\documentclass[11pt]{amsart}
\usepackage{amssymb,latexsym}
\newdimen\AAdi%
\newbox\AAbo%
\def\AAk#1#2{\s_etbox\AAbo=\hbox{#2}\AAdi=\wd\AAbo\kern#1\AAdi{}}%
\def\AAr#1#2#3{\s_etbox\AAbo=\hbox{#2}\AAdi=\ht\AAbo\raise#1\AAdi\hbox{#3}}%

\font\tenmsb=msbm10 at 12pt
\font\sevenmsb=msbm7 at 8pt
\font\fivemsb=msbm5 at 6pt
\newfam\msbfam
\textfont\msbfam=\tenmsb
\scriptfont\msbfam=\sevenmsb
\scriptscriptfont\msbfam=\fivemsb
\def\Bbb#1{{\tenmsb\fam\msbfam#1}}

\def\R{\Bbb R}

\begin{document}

\newtheorem{thm}{Theorem}
\newtheorem{lem}{Lemma}
\newtheorem{cor}{Corollary}
\newtheorem{rem}{Remark}
\newtheorem{pro}{Proposition}
\newtheorem{defi}{Definition}
\newcommand{\noi}{\noindent}
\newcommand{\dis}{\displaystyle}
\newcommand{\mint}{-\!\!\!\!\!\!\int}
\newcommand{\ba}{\begin{array}}
\newcommand{\ea}{\end{array}}

\def \bx{\hspace{2.5mm}\rule{2.5mm}{2.5mm}} \def \vs{\vspace*{0.2cm}} 
\def\hs{\hspace*{0.6cm}}
\def \ds{\displaystyle}
\def \p{\partial}
\def \O{\Omega}
\def \o{\omega}
\def \b{\beta}
\def \m{\mu}
\def \l{\lambda}
\def\L{\Lambda}
\def \ul{u_\lambda}
\def \D{\Delta}
\def \d{\delta}
\def \s{\sigma}
\def \e{\varepsilon}
\def \a{\alpha}
\def \tf{\tilde{f}}
\def\cqfd{%
\mbox{ }%
\nolinebreak%
\hfill%
\rule{2mm} {2mm}%
\medbreak%
\par%
}
\def \pr {\noindent {\it Proof:} }
\def \rmk {\noindent {\it Remark} }
\def \esp {\hspace{4mm}}
\def \dsp {\hspace{2mm}}
\def \ssp {\hspace{1mm}}

\def \u{u_+^{p^*}}
\def \ui{(u_+)^{p^*+1}}
\def \ul{(u^k)_+^{p^*}}
\def \energy{\int_{\R^n}\u }
\def \sk{\s_k}
\def \mo{\mu_k}
\def\cal{\mathcal}
\def \I{{\cal I}}
\def \J{{\cal J}}
\def \K{{\cal K}}
\def \OM{\overline{M}}

\def\fk{{{\cal F}}_k}
\def\M1{{{\cal M}}_1}
\def\Fk{{\cal F}_k}
\def\Gk{{\Gamma_k^+}}
\def\n{\nabla}
\def\uuu{{\n ^2 u+du\otimes du-\frac {|\n u|^2} 2 g_0+S_{g_0}}}
\def\uuug{{\n ^2 u+du\otimes du-\frac {|\n u|^2} 2 g+S_{g}}}
\def\sku{\sk\left(\uuu\right)}
\def\qed{\cqfd}
\def\vvv{{\frac{\n ^2 v} v -\frac {|\n v|^2} {2v^2} g_0+S_{g_0}}}
\def\vvs{{\frac{\n ^2 \tilde v} {\tilde v}
 -\frac {|\n \tilde v|^2} {2\tilde v^2} g_{^n}+S_{g_{^n}}}}
\def\skv{\sk\left(\vvv\right)}
\def\tr{\hbox{tr}}
\def\pO{\partial \Omega}
\def\dist{\hbox{dist}}
\date{October 2001}
\title[A fully nonlinear conformal flow ]
{A fully nonlinear conformal flow on locally conformally flat manifolds}
\author{Pengfei Guan}
\address{Department of Mathematics\\
 McMaster University\\
Hamilton, Ont. L8S 4K1, Canada.\\
Fax: (905)522-0935 }
\email{guan@math.mcmaster.ca}
\thanks{Research of the first author was supported in part by 
NSERC Grant OGP-0046732.}
\author{Guofang Wang}
\address{Max-Planck-Institute for Mathematics in
the Sciences\\ Inselstr. 22-26, 04103 Leipzig, Germany}
\email{gwang@mis.mpg.de}
\begin{abstract}
\noindent We study a fully nonlinear flow for conformal metrics. The
long-time existence and the sequential convergence 
of flow  are established for
 locally conformally flat manifolds. As an
application, we solve the $\sk$-Yamabe problem 
for locally conformal flat manifolds when $k \neq n/2$.
\end{abstract}
\maketitle
\section{Introduction}
Let $(M, g_0)$ be a compact, connected smooth Riemannian manifold
of dimension $n \ge 3$.
Let $S_g$ denote the Schouten tensor of the metric $g$, i.e.,
\[S_g=\frac 1{n-2}\left(Ric_g-\frac {R_g}{2(n-1)}\cdot g\right),\]
where $Ric_g$ and $R_g$ are the Ricci tensor and scalar curvature of
$g$ respectively. We are interested in deforming 
the metric in the conformal class $[g_0]$ of $g_0$ to certain extremal metric
along some curvature flow.

To introduce the flow,  let {\it $\sk$-scalar
 curvature} of $g$ (see \cite{Via1}) be
\[\sk(g):=\sk (g^{-1}\cdot S_g),\]
where $g^{-1}\cdot S_g$ is defined locally by $(g^{-1}\cdot S_g)^i_j=
g^{ik}(S_g)_{kj}$, and $\s_k(A)$ is the $k$th elementary symmetric 
function of the eigenvalues of $n\times n$ matrix $A$. 
When $k=1$, $\s_1$-scalar curvature is
just the scalar  curvature  $R$ (upto a constant multiple).
For $k>1$, it is natural to consider $\sk$ in the positive cone \[\Gamma_k^+=\{\Lambda=(\l_1,\l_2,\cdots, \l_n)\in \R^n\,|\,
\sigma_j(\Lambda)>0, \forall j\le k\}.\]
 A metric $g$ is said to be in $\Gamma_k^+$
if $\sk(g)(x) \in \Gamma_k^+, \quad \forall x\in M$. 

We propose the following curvature flow
\begin{equation}\label{flow1}
\left\{\ba{rcl}
\ds\vs \frac {d}{dt}g
&=&-(\log\sk(g)-\log r_k(g))\cdot g, \\
g(0)&=&g_0,\ea\right.
\end{equation}
where $r_k(g)$ is given by
\[ r_k(g)=  \exp \left(\frac 1{vol(g)} \int_M \log \sk (g)
\, dvol(g)\right).\]

The main goal of this paper is
to prove the following existence and convergence result of 
flow (\ref{flow1}) on locally conformally flat manifolds. 

\begin{thm}\label{thm2}
Suppose $(M,g_0)$ be a compact, connected and locally conformally flat
manifold. Assume that $g_0\in \Gamma_k^+$ and smooth, then flow
(\ref{flow1}) exists for all time $0<t <\infty$ and 
$g(t) \in C^{\infty}(M) \quad \forall t$. There exist 
two positive constants $C, \beta$ depending only on $g_0$, 
$k$ and $n$ (independent of $t$), 
such that,
\begin{equation}\label{est}
\|g\|_{C^2(M)} \le C, \quad \text{and $\forall k \neq n/2$} \quad 
\lim_{t \to \infty}\|\sk(g)-\beta\|_{L^2(M)} = 0,
\end{equation}
where the norms are taken with respect to the background metric $g_0$.
Furthermore, $\forall k \neq n/2$,
for any sequence $t_n\to \infty$, there is a 
subsequence $\{t_{n_l}\}$ with $g(t_{n_l})$ converging in $C^{1, \alpha}$-norm
($\forall 0<\alpha <1$) to some smooth 
metric $g$ and
\begin{equation}\label{eq1}\sk(g)=\beta.
\end{equation} 
\end{thm}

\medskip

If $k=1$, flow (\ref{flow1}) is a logarithmic version of the Yamabe flow
introduced by Hamilton \cite{H} in connection to the Yamabe 
problem \cite{Yamabe}. The
Yamabe problem is a problem of deforming scalar curvature to a constant 
along the conformal class, the final solution was obtained by Schoen
\cite{Schoen1} after works of Aubin \cite{Au}, Trudinger \cite{Tru} and
others. For the study of the Yamabe flow, we refer to \cite{Ye} and \cite{Chow}.

The problem of finding constant $\sk$-curvature in the conformal
class for $k \ge 2$ is a fully nonlinear version of the Yamabe problem. 
This problem was
considered in \cite{Via1}. It's reduced to a nonlinear elliptic
equation. The main difficulty in this elliptic approach is the lack of
compactness. One way to attack the problem is the blow-up 
analysis as in \cite{Schoen3} to rule out the
standard sphere. 
In the case $k=n$, Viaclovsky obtained
a sufficient condition for the solution of this nonlinear problem
in \cite{Via2}.  Chang, Gursky and Yang \cite{CGY} solved
this problem for $k=2, n=4$ on general 4 dimensional manifolds 
using a priori estimates and blow-up analysis. To do the blow-up
analysis, one needs local estimates and classification of resulting
entire solutions of the corresponding equation on $\R^n$. These type of 
local estimates
for the fully nonlinear version of Yamabe problem has been established by
us in \cite{guanwang} for general cases.
Flow (\ref{flow1}) is designed as a parabolic approach to the problem. 
As a consequence of our main 
result, we obtain the solution for the fully nonlinear version of
Yamabe problem for locally conformally flat manifolds when $k\neq \frac n2$
(of course, the $k=1$ is a result of Schoen \cite{Schoen1}).

\begin{cor}\label{yama}
If $(M,g_0)$ is a compact, connected and locally conformally flat
manifold and $g_0\in \Gamma_k^+$, $k\neq \frac n2$, then there is a 
smooth metric $g\in \Gamma_k^+$ in $[g_0]$ such that the $\sk$-curvature
of $g$ is a positive constant.
\end{cor}

\begin{rem} If $(M,g_0)$ is a compact, connected and locally conformally flat
manifold with positive Ricci curvature, then $(M, g_0)$ is conformally
equivalent to a space form. Hence, Corollary 1 is trivial for the case $k=n$.
\end{rem}

\section{Preliminaries}
In this section, we collect and derive some properties of $\sk$ and its
functionals. 

Let $(\lambda_1, \dots, \lambda_n) \in \mathbf{R}^n$.
The $k$-th elementary symmetric function is defined as 
\[\sk(\lambda_1, \dots, \lambda_n) = \sum_{i_1 < \dots < i_k}
\lambda_{i_1} \cdots \lambda_{i_k}.\]
A real symmetric $n \times n$ matrix $A$ is said to lie  in  $\Gamma_k^{+}$
if its eigenvalues lie in $\Gamma_k^{+}$. 

Let $A_{ij}$ be the $\{ij\}$-entry of an
$n \times n$ matrix. Then for $0 \leq k \leq n$, the 
$k$th {\it Newton transformation} associated with $A$ is defined to be 
\[T_k(A) = \sigma_k(A) I - \sigma_{k-1}(A) A + \cdots + (-1)^k A^k.
\]
We have (see \cite{Reilly})
\[T_k(A) ^i_j= \frac{1}{k!} \delta^{i_1 \dots i_k i}_{j_1 \dots j_k j} 
A_{i_1 j_1} \cdots A_{i_k j_k},\]
where $ \delta^{i_1 \dots i_k i}_{j_1 \dots j_k j}$ is the 
generalized Kronecker delta symbol.  Here  we use the summation convention.
By definition,
\[\sk(A)=\frac{1}{k!} \delta^{i_1 \dots i_k }_{j_1 \dots j_k } 
A_{i_1 j_1} \cdots A_{i_k j_k},\] 
\[T_{k-1}(A) ^i_j= \frac{\partial \sk(A)}{\partial A_{ij}}.\]

The following properties of $\sk$ and $\Gamma_k^+$ are well-known (e.g.,
see  \cite{Garding} ).
\begin{pro}
\label{pro1} We have
\begin{itemize}
\item [1.] Each set $\Gamma_k^{+}$ is an open convex cone.
\item [2.] $T_{k-1}(A)$ is positive definite when $A \in \Gamma_k^+$.
 \item [3.] $\log\sk$ and $\sk^{1/k}$ are concave. 
\end{itemize}
\end{pro}

The following are some variational characterizations of $\sk(g)$ (see \cite{Via1}).
\begin{pro} If (M,g) is locally conformally flat, then $T_k(S_g)$ is
divergence free with respect to the metric $g$, i.e., for any
orthonormal frame, $\forall j$
\begin{equation}\label{adj}
\sum_i \n_i (T_k(S_g)^i_j) =0.
\end{equation}
If $k \neq \frac n 2$, any solution of equation (\ref{eq1}) is a critical point of the functional
\begin{equation}\label{eq3}
{ \cal F}_k(g)=vol(g)^{-\frac{n-2k} n}\int_M \sk(g)\, dvol(g)\end{equation}
in $[g_0]$. \end{pro}

Now we consider some properties of flow  (\ref{flow1}).
\begin{lem}\label{presv}
The flow (\ref{flow1}) preserves the volume.
When $k\neq \frac n 2$ and $g_0$ is locally conformally flat,
then 
\begin{equation}\label{func}
\frac d{dt}\Fk(g) = -\frac {n-2k} 2 
vol(g)^{\frac{2k-n}2}\int_M(\sk(g)-r_k(g)) (\log \sk(g)-\log r_k(g)).
\end{equation}
Therefore, for $k <n/2$, $\frac d {dt} \fk(g(t))\le 0$, 
and for $k>n/2$, $\frac d {dt} \fk(g(t))\ge 0$.
\end{lem}

\pr 
The volume is preserved, as 
\[ \ba{rcl}
\frac{d}{dt} vol(g) & = &
\ds\vs \int_Mg^{-1} \frac d {dt} g\, dvol(g)\\
& =& \ds \int_M (\log \sk(g)-\log r_k(g)) \, dvol(g)=0.\ea\] 

On any locally conformally flat manifold, from the computation in \cite{Via1}, 
\[
\frac d{dt} \int_M \sk(g)\, dvol(g)=
\ds \frac {n-2k} 2 \int_M \sk(g) g^{-1}\cdot \frac d {dt} g \, dvol(g).\]
Since flow  preserves the
volume, we have
\[\ba{rcl}
\ds\vs \frac d{dt}\Fk(g) & = & \ds \frac {n-2k} 2 
vol(g)^{\frac{2k-n}2}\int_M \sk(g) g^{-1}
\cdot \frac d {dt} g \, dvol(g)\\
& = & \ds \frac {n-2k} 2 
vol(g)^{\frac{2k-n}2}\int_M (\sk(g)-r_k(g)) g^{-1}
\cdot \frac d {dt} g \, dvol(g)\\
\ds\vs & = & \ds -\frac {n-2k} 2 
vol(g)^{\frac{2k-n}2}\int_M(\sk(g)-r_k(g)) (\log \sk(g)-\log r_k(g)) dvol(g).
\ea \] \cqfd

If $g=e^{-2u}\cdot g_0$, one may compute that (see \cite{Via1})
\[\sk(g)=e^{2ku}\sk\left(\uuu\right).\]
Equation (\ref{flow1}) can be written in the 
following form
\begin{equation}\label{flow2}
\left\{\ba{rcl}    
\ds \vs 2\frac{\ds du}{\ds dt} & = &  \ds \log \sk
\left(\uuu\right)+2ku -\log r_k\\
u(0) & = & u_0.
\ea\right.\end{equation}

If $g=v^{-2} g_0$,
equation (\ref{flow1}) is equivalent to
\begin{equation}\label{flow3}
\left\{\ba{rcl}    
\ds \vs 2\frac{\ds dv}{\ds dt} & = &  \ds  \log \sk
\left(\vvv\right)+2k \log v-\log r_k\\
v(0) & = & v_0.
\ea\right.\end{equation}
Here all covariant derivatives are taken with respect to the fixed metric
$g_0$.

Since $g_0 \in \Gamma_k^+$, flow  (\ref{flow1}) is parabolic near
$t=0$, by the standard implicit function theorem we have 
the following short-time existence result.
\begin{pro}\label{pro2} For any $g_0\in C^2(M)$ with $\s_k(g_0) \in 
\Gamma^+_k$, there exists a positive constant $T^*$ such that
 flow (\ref{flow1}) exists and parabolic for $t \in [0, T^*)$,  
and $\forall T<T^*,$
\[ g \in C^{3,\alpha}([0,T] \times M), \forall 0<\alpha <1,
 \quad \text{and} \quad \s_k(g(t))\in \Gamma^+_k .\]
 \end{pro}

\section{A priori estimates}

In this section, we establish some a priori estimates for flow  
(\ref{flow1}). Here the assumption of the locally conformally flatness
is used. By the fundamental result of Schoen-Yau 
on positive mass theorem, if $(M,g_0)$ is not conformally equivalent to
the standard sphere, there is a precise asymptotic property for
the developing mapping of $(M,g_0)$ into the standard sphere \cite{SY}. 
Schoen applied this property together with Alexandrov's moving plane 
method developed in
\cite{GNN} to obtain the compactness in \cite{Schoen3}. Similar approach
was also taken by Ye in \cite{Ye} to obtain a Harnack type inequality
for the Yamabe flow. With little modification, 
the proof in \cite{Ye} yields the following Harnack type estimates for flow 
(\ref{flow3}).

\begin{pro}\label{pro3}
Let $(M, g_0)$ be a locally conformally flat manifold and $g_0\in \Gamma_k^+$.
If $g(t)=v^{-2}g_0$ is a solution of (\ref{flow3}) with initial metric 
$g_0$  with $\s_k(g(t)) \in \Gamma^+_k$ in a time interval $[0, T^*)$, then
\begin{equation}\label{3.1}
\sup_{M}\frac{|\n _{g_{0}} v|}{v} \le C, \quad \hbox{ for } x\in M, 
t\in [0,T^*),\end{equation}
where $C>0$ is a constant depending only on $g_0$, $k$ and $n$. 
\end{pro}
\pr
One only needs to consider the case $(M, g_0)$ is not conformally 
covered by $^n$. By \cite{SY}, there is a conformal
diffeomorphism $\Phi$ from the universal cover $\widetilde M$ of $M$ onto
a dense domain $\Omega$ of $^n$. The  boundary $\pO$ of $\Omega$ is non-empty.
Let $\pi: \widetilde M\to M$ be the covering map and $
\tilde g= (\Phi^{-1})^* \pi ^* g$ the pull-back metric. Set
$\tilde g=\tilde v^{-2}g_{^n}$, where $g_{^n}$ is the standard
metric of the unit sphere. By the conformal invariance, we have
\begin{equation}\label{3.2}
\left\{\ba{rcl}    
\ds \vs 2\frac{\ds d\tilde v}{\ds dt} & = &  \ds \tilde v^{2}\sk^{1/k}
\left(\vvs\right)-r^{1/k}_k, \quad \hbox{ in }
\Omega,\\
\tilde v(0) & = & \tilde v_0.
\ea\right.\end{equation}
Here $\n \tilde v$ and $\n^2 \tilde v$ are derivatives with respect to the 
standard metric $g_{^n}$. 
A result of \cite{SY} 
implies 
 
{\bf Asymptotic fact:}
$\forall T<T^*,
\tilde v^{-1}(x,t) \to \infty, \quad \text{uniformly $\forall t\le T$,}
\quad \hbox{ as } x\to \pO$.
\medskip

Ye used the moving plane method of \cite{GNN} to prove the Harnack
inequality for solutions of Yamabe flow. The main ingredients 
in his proof are: 
\begin{itemize}
\item [1.] The {\bf aymptotic fact};
\item [2.] Certain growth conditions at $\infty$ of the transformed function $w$ of $\tilde v$ under the stereographic projection from $^n$ to $\R^n$;
\item [2.] The invariance of the resulting parabolic equation for $w$ in $\R^n$
under reflections and translations.
\end{itemize}

In fact, the explicit form of the equation for Yamabe flow were not used
in the proof of the Harnack inequality in \cite{Ye}.
As all these ingredients are available for flow (\ref{flow3}), 
the proof in \cite{Ye} can be
adapted to get the Harnack type inequality (\ref{3.1}) for flow 
(\ref{flow3}). We will not repeat it here. \qed

By Lemma \ref{presv}, the volume is preserved under flow (\ref{flow1}). 
For $u=\log v$ satisfying (\ref{flow2}), we have 
\begin{cor} \label{cor1} There exists a  constant $C $
depending only on $(M, g_0)$, $k$ and $n$, such that for any
solution $u$ of flow  (\ref{flow2}),
\begin{equation}
\label{b1} \|u\|_{C^1} \le C.\end{equation}
\end{cor}

\medskip

Now, we consider the $C^2$-estimation for flow  (\ref{flow1}).
We note that in the case of the Yamabe flow, $C^2$ and higher 
regularity estimates follow directly from $C^1$ estimates because
the leading term is the heat equation. 
When $k\ge 2$, flow (\ref{flow1})
is fully nonlinear. It is not elliptic until the establishment of $C^2$
estimates and the positivity of $\sk(g)$ in the time interval considered.
$C^2$ estimates are crucial to the global existence of flow  (\ref{flow1}).


\begin{pro}\label{pro4}
Let $g=e^{-2u}g_0$ be a solution of flow (\ref{flow1})  
with $\s_k(g(t)) \in \Gamma^+_k$ on $M\times [0,T^*)$. Then there
is a constant $c>0$ depending only on $g_0$, $k$ and $n$ (independent of $T^*$) such that
\begin{equation}\label{b20}
|\n^2 g(t)|\le c, \quad \forall t \in [0, T^*).\end{equation}
\end{pro}

\pr For any local frame $e_1,...,e_n$, we denote $u_i=\n _i u$ and
$u_{li}=\n _l \n_i u$  the first and second covariant derivatives with 
respect to the background metric $g_0$. The similar notation will be used
for the higher order covariant derivatives.

We want to bound $\D u$. If $k\ge 2$ and $\sk(\uuu) \in \Gamma^+_k$, we know that
\[|u_{ij}|\le c_1 (\D u+ |\n u|+1).\]
Hence we only need to get an upper bound of $\D u$.
Consider $G:=(\D u +m|\n u|^2)$ on $M \times [0, T]$ 
for a given $T\in (0, T^*)$. Here $m$ is a constant to be fixed later. 
Assume $G$ achieves the maximum at $(x_0, t_0)\in  M\times [0, T]$.
Without loss of generality, we may assume that $G(x_0, t_0)\ge 1$.
Since $|\n u|$ is bounded, we may also assume that at the point 
\[G \geq \frac{1}{2}\tr(\uuu).\]
By the Newton-MacLaurin inequality, at this point
\begin{equation}\label{mean}
G \geq \frac{n}{2}(\frac{k!(n-k)!}{n!})^{\frac 1k}\sk^{\frac 1k}(\uuu).
\end{equation}

 We have, at $(x_0, t_0)$,
\begin{equation}\label{b2}
 G_t  = \sum_{l} (u_{llt}+2m u_{lt}u_l)
\ge 0\end{equation}
and
\begin{equation}\label{b2.1}
G_{i} = \sum_{l} (u_{lli}+2m u_{li}u_l) = 0.
\end{equation}
In what follows, we indicate $c$ to be the constant (which may vary line
form line) depending only on the quantities specified in the
proposition. By (\ref{b2.1}) we have at $(x_0, t_0)$,
\begin{equation}
\label{b5}
|\sum _l u_{lli} |\le c G, \quad \hbox{ for all }i.\end{equation}
Furthermore the matrix  $(G_{ij})$ is semi-negative definite
at this point. Let $W=(W_{ij})$ be a matrix defined by
\[W_{ij}=u_{ij}+u_iu_j-\frac 12|\n u|^2\d_{ij}+S(g_0)_{ij}\]
and 
$F=\log \sigma_{k}(W)$. Set
\[F^{ij}=\frac{\partial F}{\partial W_{ij}}.\]
We may assume that $ (W_{ij})$ is diagonal at $x_0$,
so $(F^{ij})$ is also diagonal at the point.
From the positivity of $(F^{ij})$, we have
\begin{equation}\label{b3}\begin{array}{rcl}
0\ge \vs \ds \sum_{i,j} F^{ij}G_{ij}
 & =&  \ds \sum_{i,j,l} F^{ij}(u_{llij} +2mu_{li}u_{lj}
 +2u_{lij}u_l).\\
\end{array}\end{equation}
Since $|\n u|$ is bounded, commutators of the covariant derivatives
can be estimated as,
\begin{equation}\label{b6}\begin{array}{rcl}
\ds\vs|u_{lij}-u_{ijl}| & \le & c,\\
\ds|u_{llij}-u_{ijll}| & \le &   c G.
\end{array}\end{equation}
In view of (\ref{b2})--(\ref{b6}) and the concavity of $F$, we have
\begin{equation}\label{b4}\begin{array}{rcl}
0&\ge & \vs\ds \sum_{i,j} F^{ij}G_{ij}\\
 &\ge & \vs \ds \sum_{i,j,l} F^{ij}(u_{ijll} +2mu_{li}u_{lj}
 +2u_{ijl}u_l)-c\sum_i F^{ii}G\\
 &=& 
 \vs \ds\sum_{i,j,l}  F^{ij}\{w_{ijll}-
 (u_iu_j-\frac 12 |\n u|^2\d_{ij}+S(g_0)_{ij})_{ll} + 2mu_{li}u_{lj}\\
 &&\vs\ds + 2mw_{ijl}u_l
 -2mu_l(u_iu_j-\frac 12 |\n u|^2\d_{ij}+S(g_0)_{ij})_l\} -c\sum_i F^{ii}G\\
 & \ge &  \vs\ds 
 \D F +2m \sum_l F_l u_l+\sum_{i,j,l} F^{ii} u^2_{jl}
 +2(m-1)\sum_{i,l} F^{ii} u^2_{li}-c\sum_i F^{ii}G\\
 & \ge &\ds \D F +2m \sum_lF_l u_l+\frac 1n G^2\sum_i F^{ii}+2(m-1)
 \sum_{i,l} F^{ii} u^2_{li}-c\sum_i F^{ii}G.
 \end{array}\end{equation} 
Recall that $ F=u_t-2ku-\log r(g)$,
(\ref{b2}) and (\ref{b4}) yield
\begin{equation}\label{b4.1}
\begin{array}{rcl}
0 & \ge & \vs\ds
(\D u_t+2m \sum_l u_l u_{tl})-2k \D u -4mk|\n u|^2\\
&& \vs \ds +\frac 1n \sum_{i}F^{ii}G^2+2(m-1) \sum_i
 F^{ii}u^2_{ii}
-c \sum_i F^{ii} G\\
&\ge & \ds \{-2k G+2(m-1) \sum F^{ii}u^2_{ii}\}+\sum_{i}F^{ii}
(\frac 1n G^2-c G).
\end{array}\end{equation}
We {\it claim} that for suitable large $m>1$,
\begin{equation}\label{4b.2}
-2k G+2(m-1) \sum F^{ii}u^2_{ii}\ge -c(1+\sum_i F^{ii}G).\end{equation}
Since $W$ and $(F^{ij})$ are diagonal at the point, $\l_i=W_{ii}$ and \[ W_i=(\l_1, \l_2, \cdots, \l_{i-1},\l_{i+1}, \cdots),\]
it is easy to check that $F^{ii}=\frac{\sigma_{k-1}(W_i)}{\s_k(W)}$.
From (\ref{b4.1})  and the identity
\[\sum_i \sigma_{k-1}(W_i)= (n-k+1)\sigma_{k-1}(W),\]
 if the {\it claim} is true,  
\begin{equation}\label{bd}
\frac{\sigma_{k-1}(W)}{\s_k(W)}(G^2-c G)\le c.
\end{equation}
The Newton-MacLaurin inequality yields 
\[\frac{\sigma_{k-1}(W)}{\s_k(W)}\ge \frac{kn}{(n-k+1)tr(\uuu)}.\]
Then the Proposition 
follows from (\ref{mean}) and (\ref{bd}).

Now we prove the {\it claim}.
From (\ref{b1}) and the assumption that
$G(x_0, t_0)\ge 1$, there is a constant $c>0$ independent of $T$ such that
\[\sum F^{ii}u^2_{ii} \ge \sum F^{ii}w^2_{ii}-c \sum F^{ii} G.\]
Since $W$ is diagonal at the point, we have the identity (e.g., see \cite{HS})
\[\sum_{i} \sigma_{k-1}(W_i)w_{ii}^2=\s_1(W)\s_k(W)-(k+1)\s_{k+1}(W).\]
In turn,
\[
\begin{array}{rcl}
\ds \sum F^{ii}w^2_{ii} & = & \vs\ds \frac 1 {\s_k(W)}\sum_{i}  
\sigma_{k-1}(W_i)w_{ii}^2\\
&=& \ds \frac 1 {\s_k(W)}(\s_1(W)\s_k(W)-(k+1)\s_{k+1}(W)).\end{array}\]
When $\s_{k+1}(W) \le 0$, the claim is automatically true. Hence
we assume that $\s_{k+1}(W)>0$. The
Newton-MacLaurin inequality
\[(k+1)\s_{k+1} \le \frac {n-k} n \s_1\s_k\]
yields 
\[\sum F^{ii}w^2_{ii} \ge \frac k n \s_1(W).\]
Now  the claim is verified if we choose
$m > n+1$.
This completes the proof of the Proposition. \qed


\section{A nonlinear eigenvalue problem}
In this section, we consider a nonlinear eigenvalue problem
\begin{equation}\label{4.1}
\sk^{1/k}\left (\n ^2 \phi + d \phi \otimes d\phi -\frac
{|\n \phi|^2} 2 +S_{g_0}\right)=\l,\end{equation}
where covariant derivatives are taken with respect to the background
metric. We say $\phi$ is  admissible 
if  $e^{-2\phi}g_0 \in \Gamma^+_k$.
We refer to \cite{Lions} and 
\cite{Wangxj} for the treatment of other types of nonlinear eigenvalue 
problems. 

\begin{thm}\label{thm4.1} Let $g_0\in \Gamma^+_k$. Then there exists a
function $\phi$ and a positive number $\l$ such that  $\phi$ is an
admissible solution of equation (\ref{4.1}). $(\phi,\l)$ is unique, 
in the sense that if there are two admissible $(\phi, \l)$
and $(\phi', \l')$ satisfying (\ref{4.1}),
then $\l=\l'$ and $\phi=\phi'+c$ for some constant $c$.
\end{thm}

In order to prove Theorem \ref{thm4.1}, we introduce an auxiliary equation
\begin{equation}\label{4.3}
\sk^{1/k}\left (\n ^2 u + d u \otimes d u -\frac
{|\n u|^2} 2 +S_{g_0}\right)=he^u+f,\end{equation}
for some functions $f>0$, $h\ge 0$. 

\begin{pro}\label{pro4.1}
Let $g_0\in \Gamma^+_k$, suppose $\frac 1L\le f \le L$ for some 
constant $L>0$ and $h\ge 0$. If $u$ is a $C^4$
admissible solution of (\ref{4.3}) and $\max_M u =\gamma$
for some constant $\gamma$, then for each $0<\alpha<1$, $l\ge 2$ integer, there is 
a constant $C$ depending only on $l$, $\alpha$, $g_0$, $L$, $\gamma$, $\|h\|_{C^l}$,
$\|f\|_{C^l}$ such that
$\|u\|_{C^{l+1,\alpha}(M)}\le C$.
\end{pro}
\noindent{\it Proof of Proposition \ref{pro4.1}.}
We will repeat the arguments in \cite{guanwang} for the a priori estimates. 
In fact, stronger local estimates hold for the solutions of (\ref{4.3})
following the same lines of the proofs in \cite{guanwang}. Here we only
concentrate on the global estimates.

We first obtain a $C^1$ bound. Since $\max_M u =\gamma$, we only need to
bound the gradient of $u$. 
Let
$W=(\n ^2 u+du\otimes du-\frac {|\n u|^2}2 g_0+S_{g_0})$ and
let $w_{ij}$ be the entries of $W$.
 Set $H=|\n u |^2$ and assume that $H$ achieves its maximum at
$x_0$. 
After appropriate choice of
the normal coordinates at $x_0$,
we may assume that
$W$ is
diagonal at the point.
Since $x_0$ is the maximum point of $H$, we have $H_i(x_0)=0$, i.e.,
\begin{equation}
\label{10}
 \sum_{l=1}^n u_{il}u_l=0.\end{equation}
We may assume that $H(x_0)\ge A^2_0$ and 
$|S_{g_0}|\le A_0^{-1}|\n u|^2$ for some 
large fixed number $A_0$ to be chosen later.
  
Since $x_0$ is the maximum point of $H$, 
the matrix 
$ ( H_{ij})=\left(2u_{lij}u_l+2u_{il}u_{jl}\right)$
is semi-negative definite.
Set 
\[ F^{ij}=\frac{\partial {\s_k}^{1/k}}{\partial w_{ij}}. \]
$(F^{ij})$ is a diagonal matrix at $x_0$ as $W$ is diagonal. 

Again, as in the proof of Proposition \ref{pro4},
we denote $C$ (which may vary from line to line)
as a constant depending only on the quantities mentioned in this proposition. 
Since $F=(he^u+f)^{1/k}$, $\sum_l F_l u_l\le C(H+1)$.  We have  
\begin{equation}\label{13}
 0\ge F^{ij}H_{ij}=F^{ij}
(2 u_{lij}u_l+2 u_{il}u_{jl}).\end{equation}
Using (\ref{10}), after commuting the covariant derivatives,
the first term in  (\ref{13}) can be estimated as follows, 
\begin{equation}\label{14}\ba{rcl}
\ds\vs \sum _{i,j,l}F^{ij}u_{ijl}u_l & \ge
& \ds \sum _{i,j,l}F^{ij}u_{ijl}u_l-C|\n u |^2 \sum_{i}F^{ii}\\
& =& \vs\ds \sum_{i,j,l}F^{ij}
\{ w_{ijl}u_l-(u_iu_j-\frac{|\n u|^2} 2\d_{ij})_l u_l\}-C|\n u |^2 \sum_{i}F^{ii}\\
\ds\vs & = & \ds  \sum_{l}F_lu_l-2\sum_{i,j,l}
F^{ij}u_{il} u_j u_l +\sum_{i,k,l} F^{ii}u_{kl}u_ku_l
-C|\n u |^2 \sum_{i}F^{ii}\\
 &\ge& \ds -C|\n u |^2 \sum_{i}F^{ii}.\\
\ea\end{equation}
Here we have used the homogeneity $\sum_i F^{ii}w_{ii}=F$.
 
By Lemma 1 in \cite{guanwang}, for $A_0$ sufficiently large (depending
only on $k$, $n$, and $\|g_0\|_{C^3}$)
\begin{equation}\label{15}\sum_{i,j,l}F^{ij}u_{il}u_{jl} 
\ge A_0^{-{\frac 34}}|\n u|^4\sum_{i\ge 1}  F^{ii}.
\end{equation} 
(\ref{14}) and (\ref{15}), together with (\ref{13}) yield the desired
$C^1$ estimate.

With the $C^1$ bound, a $C^2$ bound can be obtained following the same
lines of proof of Proposition \ref{pro4} in the previous section. We
may assume $k\ge 2$, since $C^2$ bound for $k=1$ follows from 
linear elliptic theory. Since $u$ is  admissible, in this case
we only need to get an upper bound for $\D u$ as in the proof of
Proposition \ref{pro4}.
We estimate the maximum of $G= (\D u+|\n u|^2)$. We note that
$F=\s_k^{1/k}(\uuu)$ is also concave.
At any maximum point $y_0\in M$, by the similar computations as in (\ref{b4})
and (\ref{b4.1}) in the  proof of Proposition \ref{pro4}, we have,

 \begin{equation}\ba{rcl} \label{33}
0 \ge \D F + 2\sum_lF_lu_l+ \sum_{i,k,l} F^{ii} u_{kl}^2
  -C(1+G)\sum_i F^{ii}.
 \ea\end{equation}
We estimate the terms on the right hand side.
As $F=(he^u+f)^{1/k}$, we have
\[\sum_lF_lu_l \ge -C, \quad \sum_{l\ge 1} F_{ll} \ge -CG .\]
By the facts $\sum F^{ii}\ge 1$ for $F=\sk^{\frac 1k}$ and $G(y_0)\ge 1$, 
the above yields 
\begin{equation}\label{34}\ba{rcl}
0&\ge& 
\ds\vs 
-C \sum F^{ii} G + \sum_{i,k,l} F^{ii} u_{kl}^2 \ge
- C\sum F^{ii}G
+ \frac 1n \sum_{i} F^{ii} (\D  u)^2 \\
 &\ge& 
\ds\vs \sum F^{ii}\left\{-CG+\frac 1{n} G^2\right\}.
\ea\end{equation}
It follows from (\ref{34}) that at $y_0$,
$G\le C$. 

The higher regularity estimates follow from the 
Evans-Krylov theorem (e.g., \cite{Krylov}). \qed

\noindent{\it Proof of Theorem \ref{thm4.1}.}
First we want to prove that for small $\l>0$ the following equation
has a unique smooth admissible solution
\begin{equation}\label{4.43}
\tilde F(u)=:\sk^{1/k}\left (\n ^2 u + d u \otimes d u -\frac
{|\n u|^2} 2 +S_{g_0}\right)-e^u=\l.\end{equation}
Since $\frac {\partial \tilde F}{\partial u}<0$, the uniqueness for the solutions of (\ref{4.43})
follows from the Maximum principle. For the same reason, the kernel 
of the linearized operator of $\tilde F$ is trivial at any admissible $v$.
We note that for any admissible $u$ and $\forall 0\le t \le 1$, 
since $g_0 \in \Gamma_k^+$, 
$u_t=tu$ is also admissible. The linearized operator 
of $\tilde F$ at $u_0=0$ is 
$L(\rho)=\tr \{ T_{k-1}(S_{g_0}) \n^2_{g_0}\rho\} -h \rho$.
By (\ref{adj}), it is
self-adjoint 
with respect to the metric $g_0$. Since the index of elliptic 
operator is invariant 
under homotopy, so the kernel of
the adjoint operator of the linearized operator of $\tilde F$ at $u$
is trivial.
That is, the linearized operator of $\tilde F$ is invertible at every
admissible solution $u$. This fact will be used later in the proof.

We now want to show the existence using the
continuity method.  
Since $g_0\in \Gamma^+_k$, there is a constant $C_1>1$
such that $C_1^{-1}<\sk^{1/k}(S_{g_0})<C_1$. Thus, for small $\l>0$,
one can find two constants $\underline \d < 0 < \overline \d$
such that
\[e^{\underline \d}+\l <\sk^{1/k}(S_{g_0})<e^{\overline \d}+\l.\]
Let $v=\underline \d$, we have $\tilde F(v)=f$ for some smooth
positive function $f\ge \l$. For $0\le t \le 1$,
let us consider the equation
\begin{equation}\label{4.44}
\tilde F(u)=t\l+(1-t)f.
\end{equation}
By the Maximum principle, for any solution $u$ of (\ref{4.44}),
$\min u \ge \min v=\underline \d$. Also for $\tilde v=\overline \d$,
as $\tilde F(\tilde v)<\l$, again by the Maximum principle, 
$\max u \ge \max \tilde v=\overline \d$. That is, $u$ is bounded.
By Proposition \ref{pro4.1} we have the uniform a priori estimates
for solution $u$ of (\ref{4.44}). This implies the closeness. The 
openness follows from the standard implicit function theorem since
the linearized operator of $\tilde F$ is invertible. Therefore, the
existence of the unique solution of (\ref{4.43}) is established for
$\l>0$ small.

Set
\[\L =
\{\l>0\,|\, (\ref{4.43}) \hbox{ has a solution}\}.\] 
Since $\L \neq \emptyset$, we define
\[\l^*=
\sup _{\l\in\L}\l.\]
We claim $\l^*$ is finite.
For any admissible solution $u$ of (\ref{4.43}), by the
 Newton-MacLaurin inequality, 
\begin{equation}\ba{rcl}\label{4.5}
\ds  \l &< & \vs \ds e^u+\l   =  \sk^{1/k}(\uuu)\\ 
 &\le& \ds \vs \frac{(\frac{n!}{k!(n-k)!})^{\frac 1k}}{n} \sigma_1(\uuu)\\
&=& \ds \vs \frac{(\frac{n!}{k!(n-k)!})^{\frac 1k}}{n}
( \D u -\frac {n-2}2 |\n u|^2) +cR_0,
\ea
\end{equation}
where $R_0$ is the scalar curvature function of $(M,g_0)$ and $c$ is a 
constant depending only on $n$ and $k$.
Integrating the above inequality over $M$, we get 
$\l \le c \frac{\int_M R_0d vol(g_0)} {Vol(g_0)}$.
We conclude that $\l^*\le c \frac{\int_M R_0d vol(g_0)} {Vol(g_0)} $.

For any sequence$\{\l_i\} \subset \Lambda$ with $\l_i \to \l^*$, and let 
$u_{\l_i}$ be the corresponding solution of (\ref{4.43}) 
with $\l=\l_i$, $i=1, 2,3,\cdots.$ 
We claim that $\max_M u_{\l_i} \to  -\infty$ as $i \to \infty$. 
Suppose $\max_M u_{\l_i} \ge -C_0$.
From equation (\ref{4.43}), at any maximum point $x_0$ of $u_{\l_i} $,
$\max_M u_{\l_i} \le -C$ for some constant $C$ depending only on $n$, $k$,
$g_0$. Then Proposition \ref{pro4.1} implies that $ u_{\l_i}$
(by taking a subsequence)
converges to a smooth function $u_0$ in $C^4$, such that $u_0$ satisfies
(\ref{4.43}) for $\l=\l_1$. 
As remarked at the beginning of the proof,
the linearized operator of  equation (\ref{4.43})
is invertible. By the standard implicit function theorem, we have a
solution of (\ref{4.43}) for $\l=\l_1+\e$ for $\e>0$ small, this is a contradiction. 
Hence $\max_M u_{\l_i}
\to -\infty$ as $i\to \infty$. Now let
$w_{\l_i} = u_{\l_i} -\max_Mu_{\l_i}.$
It is clear that $w_{\l_i}$ satisfies, 
\[\sk^{1/k}\left(\n^2 w_{\l_i}+d w_{\l_i}\otimes d w_{\l_i}
-\frac {|\n w_{\l_i}|^2} {2} g +S_g\right) = e^{\max_M u_{l_i}} e^w+\l_i\]
with $\max_Mw_{\l_i} \to 0$. By Proposition \ref{pro4.1} again, 
$w_{\l_i}$ converges to a smooth function $\phi$ in $C^4$ and $\phi$ satisfies
(\ref{4.1}) with $\l =\l^*$. 

Finally we prove the uniqueness. For each admissible $u$, let
\[
W=(\uuu)\] 
and 
\[a_{ij}(W)=\frac{\partial \sk(W)}{\partial w_{ij}}.\] 
For any smooth functions $u_0$ and $u_1$, 
let $v=u_1-u_0$, $u_t=t u_1 +(1-t)u_0$ and 
$W_t=\D u_t +du_t \otimes d u_t -\frac {|\n u_t|^2} {2} g +S_g$.
By (\ref{adj}), the
following identity holds
\begin{equation}\label{dvi}
\sk(W_1)-\sk(W_0)=\sum_{ij} \n_j [(\int_0^1 a_{ij}(W_t)dt)\n_{i}v]
+\sum_l b^l \n_l v,
\end{equation}
for some bounded functions $b^l$, $l=1,..., n$.
If $u_0=\phi$ and $u_1=\phi'$ are two admissible solutions of (\ref{4.1}) for
some $\l$ and $\l'$ respectively, then $(\int_0^1 a_{ij}(W_t)dt)$ is positive definite. Therefore, $\phi=\phi' +c$ for some constant $c$
by (\ref{adj}) and the Maximum principle. 
The proof of Theorem \ref{thm4.1} is complete.\qed

\section{Existence and convergence}

First, we want to use the a priori estimates established in the
section 3 to obtain the long time
existence of flow . These estimates are independent of time
$t$, as long as $u(t,x)$ stays in the positive cone $\Gamma^+_k$. 
Therefore, to establish the long time existence, we need to show
$\s_k$ stays strictly positive for all $t$.

\begin{lem}\label{lem4} 
We have
\begin{equation} \label{7}
\frac d {dt} \log \sk (g) =\ds \frac 1{2\s_k(g)}\tr \{ T_{k-1}(S_g) 
\n^2_g \log \sk(g)\} +\log\sk(g)-\log r_k(g).
\end{equation}
\end{lem}

\pr
It is easy to check, see for example \cite{Reilly},
\[\ba{rcl}
\ds \frac d{dt} \sk(g)  
 & =  & \ds k\sk(g)g \cdot \frac  d {dt} (g^{-1})
 + \tr \{ T_{k-1}(S_g)
g^{-1}\frac d {dt}  S_g\} .\ea\]
Under the conformal change 
$g=e^{-2u}g_0$, the Schouten tensor is changed as follows
\[S_g  = \uuu.\]
Hence, we have
\[\frac d {dt} S_g= \n^2_g \frac {du}{dt}=-\frac 12\n^2_g (g^{-1} \frac d {dt}
g ).\]
Now we have
\[\ba{rcl}
\ds \vs \frac d{dt} \log \sk(g)& = & \ds -\frac 1{2}\sk^{-1} 
 \tr \{ T_{k-1}(S_g)
g^{-1} \n^2_g (g^{-1} \frac d {dt}
g ) \}+\log\sk(g)-\log r_k(g)\\
&=&  \ds \frac 1{2\sk(g)}\tr \{ T_{k-1}(S_g) 
\n^2_g \log \sk(g)\} +\log\sk(g)-\log r_k(g). \ea\]
\qed

\begin{cor}\label{pos}
Suppose $g(t)$ is a solution of equation (\ref{flow1}) with initial 
metric 
$g_0$ on a time interval $[0, T^*)$, then there is a positive constant 
$c>0$
depending on $g_0$, $k$, $n$ (independent of $t$), such that 
$\forall t\in [0, T^*)$, 
\begin{equation}
\sk(g(t)) \ge c e^{-\frac{e^t}{c}}.
\end{equation}
Therefore, flow  (\ref{flow1}) exists for all $t >0$.
\end{cor}

\pr
By Propositions
\ref{pro3} and \ref{pro4},  $\log r_k(g)\le C$ for some constant $C$
independent of $t$. Let $f=e^{-t} \log \sk(g)-Ce^{-t}$. By
Lemma \ref{lem4} $f$ satisfies
a parabolic inequality
\[\frac d {dt} f \ge \frac 1{2\s_k(g)}\tr \{ T_{k-1}(S_g) 
\n^2_g f\}.\]
Now the corollary follows from the maximum principle.
Again by Propositions \ref{pro3}, \ref{pro4}, and what we just proved, 
flow (\ref{flow1}) stays uniformly parabolic
in $[0, T)$ for any fixed $T< \infty$. By Krylov Theorem for fully nonlinear parabolic equations \cite{Krylov}, 
$g(t)\in C^{2, \alpha_T}, \forall 0\le t \le T$ 
for some $\alpha_T>0$. In turn, $\|g\|_{C^{l}([0,T]\times M)} \le C_{l, T}$
for all $l\ge 2$. Now the long time existence follows from 
the standard parabolic PDE theory. \qed

The next proposition will also be used in our proof of Theorem \ref{thm2}. 
\begin{pro}\label{cor4} If $g_0\in \Gamma^k_+$, then there is no
$C^{1,1}$ function $u\in \overline{\Gamma^k_+}$ 
such that $\sk(g)=0$ for $g=e^{-2u}g_0$ in the viscosity sense.
\end{pro}

\noindent{\it Proof of Proposition \ref{cor4}}.
By Theorem \ref{thm4.1}, there is a smooth admissible solution $\phi$ of
(\ref{4.1}) for some $\l>0$. Assume by contradiction that there is a $C^{1,1}$ solution $u$
of
\[\sk \left(\uuug\right)=0.\]
Let $u_0=u$ and $u_1=\phi$ and $v=u_1-u_0$, then (\ref{dvi}) still holds for
$v$ in weak sense (since $u \in C^{1,1}$). Also, 
$(\int_0^1 a_{ij}(W_t)dt)$ is uniformly positive definite everywhere 
since $ a_{ij}(W_t)$ is semi-positive
definite for all t and uniformly positive definite for $\delta \le t \le 1$,
$\forall \delta >0$. By Maximum principle again, $u=\phi +c$ for some
constant $c$. That is, 
\[0=\sk \left(\uuug\right)=\sk(W_1)=\l>0.\]
This is a
contradiction. \qed

\bigbreak

{\noindent \it Proof of Theorem \ref{thm2}}. 
{The} a priori estimates and the longtime existence  are
proved in Proposition \ref{pro4} and
Corollary \ref{pos}. We now prove the convergence 
for $k \neq n/2$. Since the volume is preserved under
flow (\ref{flow1}),
from  (\ref{func}) we have $\forall T$
\[\int_{0}^T \int_M(\sk(g)-r_k(g))(\log \sk (g)-\log r_k)dvol(g) dt 
\le \frac{2V^{\frac{n-2k}{2}}(g_0)}{n-2k}|{\cal F}_k(g(T))- 
{\cal F}_k(g_0)|.\]
By Proposition \ref{pro4}, $\sk(g)$ and ${\cal F}_k(g)$ are 
 uniformly bounded. We conclude that
\begin{equation}\label{conv}
\int_{0}^\infty\int_M(\sk(g)-r_k(g))^2dvol(g) dt  <\infty.
\end{equation}
Therefore for any sequence $t_n \to \infty$, there is a
subsequence 
$\{t_{n_l}\}$) such that $r_k(g_{t_{n_L}}) \to \beta$, for some 
nonnegative constant $\beta$, and 
$\int_M (\sk(g_{t_{n_l}})-r_k(g_{t_{n_l}}))^2dvol(g_{t_n}) \to 0$.
Again, from Proposition \ref{pro4} (by taking subsequence) 
$ u(t_{n_l})$ converges to some $C^{1,1}$ function $u_\infty$ 
in $C^{1,\alpha}$ norm
for any $0<\alpha <1$. By Lebesgue domination theorem, 
\begin{equation}\label{final}
\s_k(g_{\infty})=\beta, \quad \text {almost everywhere},
\end{equation}
where $g_{\infty}=e^{-2u_{\infty}}g_0$. Since $u_{\infty} \in C^{1,1}$,
$u_{\infty}$ satisfies the equation (\ref{final}) in viscosity sense.
By Lemma \ref{presv} $\Fk(g(t)) $
is monotonic, $\beta$ is independent of the choice of the sequence. 
It follows (\ref{conv}) again that 
\[\lim_{t \to \infty} \|\sk(g)-\beta\|_{L^2(M)} =0.\]
Finally, $\beta$ is positive by (\ref{final}) and Corollary \ref{cor4}. 
We conclude that $u_{\infty} \in C^{\infty}(M)$ by the Evans-Krylov Theorem.
We remark that the positivity of $\beta$
for $k>n/2$ also follows from Lemma \ref{presv}, Propositions \ref{pro3} and \ref{pro4}, since $r_k(g_t)$ is  uniformly bounded from below. \qed

\medskip
\begin{rem} With some modifications in our proof, we can establish
the similar a priori estimates for flow  ($0 \le l<k\le n$)
\begin{equation}\label{qflow1}
\left\{\ba{rcl}
\ds\vs \frac {d}{dt}g
&=&-(\log(\frac{\sk(g)}{\sigma_l(g)})-\log r_{k,l}(g))\cdot g, \\
g(0)&=&g_0,\ea\right.
\end{equation}
with $r_{k,l}(g)$  given by
\[ r_{k,l}(g)=  \exp \left(\frac 1{vol(g)} \int_M \log (\frac{\sk(g)}{\sigma_l(g)}) 
\, dvol(g)\right).\]
Therefore the global existence result is valid for (\ref{qflow1}).
\end{rem}
\begin{rem} When $k=\frac n2$, the right hand side in (\ref{func}) is 
trivial. We do not know the behavior of flow (\ref{flow1}) at $\infty$
in this case. The sequential convergence of flow (\ref{flow1}) is suffice for
Corollary 1.
It is of interest to know the uniqueness of $g_\infty$ and  
general convergence of flow (\ref{flow1}).
Finally, we would like to ask whether Theorem \ref{thm2} is 
ture without the assumption of locally conformally flatness.
\end{rem}


\begin{thebibliography}{99}
\bibitem{Au} T. Aubin, {\em \'Equations diff\'erentilles non lin\'eaires
et probl\'eme de Yamabe concernant la courbure scalaire},
J. Math. Pures Appl. \textbf{55} (1976), 269-296.

\bibitem{CGY} A. Chang, M. Gursky and P. Yang, {\em An a priori 
estimate for a fully nonlinear equation on four-manifolds},
 Preprint, August, 2001.
\bibitem{Chow} B. Chow, {\em The Yamabe flow on locally conformally
flat manifolds with positive Ricci curvature,} Comm. Pure Appl. Math. 
\textbf{45}
(1992), no. 8, 1003--1014.
\bibitem{Garding}
L. G{\.a}rding, \emph{An inequality for hyperbolic polynomials}, J. Math.
  Mech. \textbf{8} (1959), 957--965.
\bibitem{GNN}B. Gidas, W.-M. Ni and  L. Nirenberg,
\emph{Symmetry and related properties via the maximum principle}, 
Comm. Math. Phys. \textbf{68} (1979), no.~3, 209--243.
\bibitem{guanwang} P. Guan and G. Wang,
  Local estimates for a class of conformal equation
  arising from conformal geometry, Preprint, Aug. 2001.
\bibitem{H} R. Hamilton, {\em The Ricci flow on surfaces,} Mathematics
and general relativity, Contemporary Math., Vol. \textbf{71}, AMS, (1988), 237-262.
\bibitem{HS} G. Huisken and C. Sinestrari,
{\em Convexity estimates for mean curvature flow and singularities of mean
convex surfaces}. Acta Math. \textbf {183},
(1999), 45--70.
\bibitem{Krylov} N. Krylov, {\em Nonlinear elliptic and parabolic equations
of the second order,} D. Reidel, (1987)
\bibitem{Lions}  P. L. Lions, {\em Two remarks on the Mpnge-Amper\'e},
Ann. Mat. Pura Appl. \textbf{142} (1985), 263-275
\bibitem{Reilly} R. Reilly, {\em On the Hessian of a function and the
curvatures of its graph,} Michigan Math. J. \textbf{20} (1973), 373--383.
\bibitem{Schoen1}R. Schoen, {\em Conformal deformation of a Riemannian
metric to constant scalar curvature,} J. Differential Geom. \textbf{20} (1984), no. 2,
479--495.
\bibitem{Schoen3} R. Schoen, {\em On the number of constant
scalar curvature metrics in a conformal class,} Differential
geometry, 311--320, Pitman Monogr. Surveys Pure Appl. Math.,
\textbf{52}, Longman Sci. Tech., Harlow, 1991. 
\bibitem{SY} R. Schoen and S. T. Yau,
   {\em  Conformally flat manifolds,
Kleinian groups and scalar curvature,} Invent. Math. \textbf{92} (1988), no. 1,
47--71.
\bibitem{Tru} N. Trudinger, {\em Remarks concerning the 
conformal deformation of Riemannian structures on compact manifolds}, 
Ann. Scuola Norm. Sup. Pisa \textbf{22}, (1968), 265-274.
\bibitem{Via1} J. Viaclovsky, {\em Conformal geometry, contact geometry and the calculus of
variations}, Duke Math. J. \textbf{101} (2000), no. 2, 283--316.
\bibitem{Via2} J. Viaclovsky, {\em Estimates and some existence
results for some fully nonlinear elliptic equations on Riemannian
manifolds,} {\bf math-DG}/0104227.
\bibitem{Wangxj} 
X. J.  Wang, {\em A class of fully nonlinear elliptic equations and related
 functionals}, Indiana Univ. Math. J. \textbf{43} (1994), no.~1, 25--54.
\bibitem{Yamabe} H. Yamabe,
{\em On a deformation of Riemannian structures on compact manifolds,} 
Osaka Math. J. \textbf{12}, 1960 21--37.
\bibitem{Ye} R. Ye, {\em Global existence and convergence of
Yamabe flow,} J. Differential Geom. \textbf{39} (1994), no. 1, 35--50
\end{thebibliography}
\end{document}